\documentclass[11pt]{amsart}
\usepackage[spanish]{babel}
\usepackage[utf8]{inputenc}

\newcommand{\titulo}[1]{\begin{center}
\huge\bf{#1}
\end{center}}
\newcommand{\autor}[1]{\begin{center}
\large\bf{#1}
\end{center}\vspace*{0.3cm}}
\newcommand{\resumen}[1]{\begin{center}\large{\textbf{Resumen}}\end{center}
\begin{changemargin}{1cm}
{\textit{#1}}
\end{changemargin}\vspace{0.8em}}
\usepackage[margin=2cm]{geometry}
\newtheorem{teorema}{Teorema}
\newtheorem{lema}{Lema}[teorema]
\newtheorem{corolario}{Corolario}[teorema]

\theoremstyle{definition}

\spanishdecimal{.}
\newtheorem*{prueva}{Dem}

\usepackage{hyperref}
\usepackage[mathscr]{euscript}
\usepackage{mathrsfs}
\usepackage{yfonts}
\usepackage{ragged2e}
\usepackage{subcaption}
\usepackage{pgf, tikz}
\usepackage{enumerate}
\usepackage{ytableau} 
\usepackage[all]{xy} 
\usepackage{ bbold }

\def\tb#1{\textbf{#1}}

\def\U{\mathbb{1}}

\def\N{{\mathbb N}}
\def\D{\mathscr{D}}
\def\C{\mathscr{C}}
\def\J{\ensuremath\mathcal{J}}

\def\T{\mathscr{T}}
\def\Tx{\Hat{\mathscr{T}}}
\def\Ty{\Tilde{\mathscr{T}}}
\def\I{\mathscr{I}}
\def\H{\mathcal{H}}

\newcommand{\Ct}[1]{\tb{Cat}_{#1}}

\def\rg{\textcolor{red}}

\def\a{\alpha}
\def\l{\lambda}
\def\e{\eta}

\def\vf{{\varphi}}

\def\g{\gamma}

\def\b{\beta}

\def\de{\delta} 
\def\De{\Delta}

\newcommand{\cmb}[2]{  \binom{ #1}{#2} }

\def\Bz{\mathcal{B}}
\newcommand{\Biz}[3]{\dfrac{1}{#1+#2}\cmb{#3#1+#3#2}{#3#1}}
\newcommand{\Bizu}[3]{\Bz_{#3_{1}}^{(#1,#2)}\Bz_{#3_{2}}^{(#1,#2)}\cdots\Bz_{#3_{l}}^{(#1,#2)}      }

\def\tb#1{\textbf{#1}}

\newcommand{\dint}[1]{\left\lfloor #1 \right\rfloor}

\newcommand{\dsum}[3]{ \displaystyle\sum\limits_{#1}^{#2} { #3}  }

\begin{document}
\titulo{Caminos de Dyck contenidos\\en un diagrama de Ferrers}
\autor{Jos\'{e} Eduardo Bla\v{z}ek$^1$} 

\center{\small ${}^1$Laboratoire de Combinatoire et d'Informatique Math\'{e}matique\\Universit\'{e} du Qu\'{e}bec \`{a} Montr\'{e}al\\
\url{jeblazek@lacim.ca}
}

\

\resumen{
Se considerara el problema de contar el subconjunto de caminos de Dyck contenido en un diagrama de Ferrers. Esta enumeración atañe al conteo de los elementos en una rama del árbol de Kréwéras. Mediante el uso de diagramas de Ferrers asociados a los caminos de Dyck, hemos desarrollado métodos de comparación y de descomposición de diagramas para obtener fórmulas enumerativas en términos de números de Catalan. Estos métodos han sido desarrollados en forma de algoritmos y codificados en SAGE para su verificación.} 

\section{Introducción}
\justify{
El estudio de los caminos Dyck es un tema central en la combinatoria, ya que proporcionan una de las muchas interpretaciones de los números de Catalan. Una visión parcial se puede encontrar, por ejemplo, en la presentación integral de Stanley de la combinatoria enumerativa \cite{Stanley} (vea también \cite{BLL}). Como lenguaje generado por una gramática algebraica se caracteriza en términos de un lenguaje de Dyck, que son importantes en la informática teórica, como en \cite{Eilenberg}. En un alfabeto de dos letras que corresponden a las expresiones entre paréntesis y así se pueden interpretar en términos de caminos en un cuadrado.
Entre las muchas generalizaciones posibles, es natural considerar caminos en un rectángulo, véase, por ejemplo Labelle y Yeh \cite{LY90}, y más recientemente Duchon \cite{Duchon} o Fukukawa \cite{Fu13}. En la combinatoria algebraica los caminos Dyck están relacionados con las funciones de Parking y la teoría de la representación del grupo simétrico \cite{GMV14}. En otras ciencias como en la biología, los encontramos juntos a las particiones no cruzadas y a los árboles explicando las estructuras del ARN, por ejemplo en los trabajos Schmitt y Waterman \cite{SW94},  Simion \cite{Sim00}, o Bakhtin y Heitsch\cite{BH09}.  La motivación para el estudio de estos objetos se debe a un intento de comprender mejor los vínculos entre ellos. En otras palabras, en el presente trabajo mostramos un método para obtener formulas enumerativas para los caminos situados en ramas del árbol de Kréwéras \cite{Krew65} usando los números de Catalan y los diagramas de Ferrers. Mas precisamente, los resultados principales de este articulo son el {\bf método de comparación de diagramas}, el {\bf método de descomposición} y el {\bf desarrollo de formulas enumerativas} en términos de números de Catalan.
}

\section{Definiciones y notaciones}

Hemos utilisado la notación de Lotario \cite{Lothaire1}. Una \emph{alfabeto} es un conjunto finito $ \Sigma$, cuyos elementos se llaman \emph{letras}. El conjunto de palabras finitas más
$\Sigma$ se denota $\Sigma^*$ y $\Sigma^+=\Sigma^*\setminus \{\varepsilon\}$ es el conjunto de palabras no vacíos donde $\varepsilon \in\Sigma^*$ es la palabra vacía.
El número de apariciones de una letra dada $\alpha$ en la palabra $w$ se denota $|w|_\alpha$ y
$|w| = \sum_ {\alpha \in \Sigma} | w |_\alpha $ es la longitud de la palabra.\\

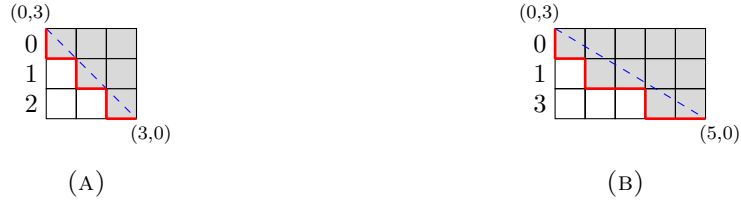
\begin{figure}[h] 
\centering
\subcaptionbox{}[0.4\linewidth]{
\begin{tikzpicture}[scale=0.4, every node/.style={scale=0.9}]
 \draw ( 0, 0 ) node{0}; 
\draw ( 0, -1 ) node{1};
\draw ( 0, -2 ) node{2};  
\draw ( 0,1 ) node{$_{(0,3)}$}; 
\draw ( 4,-3 ) node{$_{(3,0)}$};  
 \draw  ( 0.50,-2.5 ) rectangle ( 1.5, -1.5 );
\draw  ( 1.5, -2.5 ) rectangle ( 2.5, -1.5 ) ; 
\draw  [fill=gray!30 ]( 2.5, -2.5 ) rectangle (3.5, -1.5 ); 
\draw  ( 0.50, -1.5 ) rectangle ( 1.5, -0.50 ); 
\draw [fill=gray!30 ] ( 1.5, -1.5 ) rectangle (2.5, -0.50 ); 
\draw [fill= gray!30 ] (2.5, -1.5 ) rectangle ( 3.5,-0.50 );
\draw [fill=gray!30 ] ( 0.50, -0.50 ) rectangle (1.5, 0.50 ); 
\draw [fill= gray!30 ] (1.5, -0.50 ) rectangle ( 2.5,0.50 ); 
\draw [fill= gray!30 ] ( 2.5,-0.50 ) rectangle ( 3.5, 0.50 );
 \draw ( 1.5, 0.50 ) -- (1.5, -2.5 ); 
 \draw ( 2.5,0.50 ) -- ( 2.5, -2.5 ); 
 \draw (3.5, 0.50 ) -- ( 3.5,-2.5 ); 
 \draw [color= blue, dashed]+(0.5,0.5) -- ( 3.5, -2.5 ); \draw [line
width=1pt,color=red] ( 0.50, 0.50 ) -- (0.50, -0.50 ); \draw [line
width=1pt,color=red] ( 0.50, -0.50 ) -- (1.5, -0.50 ); \draw [line
width=1pt,color=red] ( 1.5, -0.50 ) -- (1.5, -1.5 ); \draw [line width=1pt,color=red]( 1.5, -1.5 ) -- ( 2.5,-1.5 ); 
\draw [line width=1pt,color=red] ( 2.5, -1.5 ) -- ( 2.5, -2.5 ); 
\draw[line width=1pt,color=red] ( 2.5, -2.5 ) -- (3.5, -2.5 ); 
\end{tikzpicture}
}
\subcaptionbox{}[0.4\linewidth]{
\begin{tikzpicture}[scale=0.4, every node/.style={scale=0.9}]
 \draw ( 0, 0 ) node{0}; 
\draw ( 0, -1 ) node{1};
\draw ( 0, -2 ) node{3}; 
\draw ( 0,1 ) node{$_{(0,3)}$}; 
\draw ( 6,-3 ) node{$_{(5,0)}$};  
\draw  ( 0.50,-2.5 ) rectangle ( 1.5, -1.5 );
\draw  ( 1.5, -2.5 ) rectangle ( 2.5, -1.5 ) ; 
\draw  ( 2.5, -2.5 ) rectangle (3.5, -1.5 ); 
\draw [fill= gray!30 ] (3.5, -2.5 ) rectangle ( 4.5,-1.5 ); 
\draw [fill= gray!30 ] ( 4.5,-2.5 ) rectangle ( 5.5, -1.5 );
\draw  ( 0.50, -1.5 ) rectangle ( 1.5, -0.50 ); 
\draw [fill=gray!30 ] ( 1.5, -1.5 ) rectangle (2.5, -0.50 ); 
\draw [fill= gray!30 ] (2.5, -1.5 ) rectangle ( 3.5,-0.50 );
 \draw [fill= gray!30 ] ( 3.5,-1.5 ) rectangle ( 4.5, -0.50 );
\draw [fill= gray!30 ] ( 4.5, -1.5 ) rectangle ( 5.5, -0.50 ); 
\draw [fill=gray!30 ] ( 0.50, -0.50 ) rectangle (1.5, 0.50 ); 
\draw [fill= gray!30 ] (1.5, -0.50 ) rectangle ( 2.5,0.50 ); 
\draw [fill= gray!30 ] ( 2.5,-0.50 ) rectangle ( 3.5, 0.50 );
\draw [fill= gray!30 ] ( 3.5, -0.50 )rectangle ( 4.5, 0.50 ); 
\draw [fill= gray!30] ( 4.5, -0.50 ) rectangle ( 5.5, 0.50 );
 \draw ( 0.50,-0.50 ) -- ( 5.5, -0.50 );
 \draw( 0.50, -1.5 ) -- ( 5.5,-1.5 ); 
 \draw ( 1.5, 0.50 ) -- (1.5, -2.5 ); 
 \draw ( 2.5,0.50 ) -- ( 2.5, -2.5 ); 
 \draw (3.5, 0.50 ) -- ( 3.5,-2.5 ); 
 \draw ( 4.5, 0.50 ) -- (4.5, -2.5 );
 \draw [color= blue, dashed]+(0.5,0.5) -- ( 5.5, -2.5 ); \draw [line
width=1pt,color=red] ( 0.50, 0.50 ) -- (0.50, -0.50 ); \draw [line
width=1pt,color=red] ( 0.50, -0.50 ) -- (
1.5, -0.50 ); \draw [line
width=1pt,color=red] ( 1.5, -0.50 ) -- (
1.5, -1.5 ); \draw [line width=1pt,color=red]
( 1.5, -1.5 ) -- ( 2.5,
-1.5 ); \draw [line width=1pt,color=red] ( 2.5
, -1.5 ) -- ( 3.5, -1.5 ); \draw
[line width=1pt,color=red] ( 3.5, -1.5 ) -- (
3.5, -2.5 ); \draw [line width=1pt,color=red]
( 3.5, -2.5 ) -- ( 4.5,
-2.5 ); \draw [line width=1pt,color=red] ( 4.5
, -2.5 ) -- ( 5.5, -2.5 );
\end{tikzpicture}
}
\caption{Caminos de Dyck y diagramas de Ferrers.}\label{E35D}
\end{figure}

\paragraph{\bf Palabras y caminos de Dyck.} Es bien conocido que el lenguaje de las palabras Dyck en  $\Sigma=\{\tb0,\tb1\}$ es un lenguaje generado por la gramática algebraica $D\rightarrow \tb0D\tb1D +\varepsilon$. Estos son enumerados por los números de Catalan (ver \cite{KOS09}),
\[\Ct{n}= \frac{1}{n+1}{2n \choose n},\]
y pueden ser interpretados como los caminos inscriptos en un retículo de forma cuadrada que no cruzan la diagonal, $n\times n$ usando pasos hacia la derecha y hacia bajo (see Fig. \ref{E35D}(A)).\\

\vspace{-0.5cm}

\

Más precisamente, un camino $(a,b)$-Dyck es un camino sudeste en un retículo, comenzando en $(0,a)$ y terminando en $(b,0)$,que se mantiene por debajo de la $(a,b)$-diagonal, siendo esta el segmento de une los puntos $(0,a)$ y $(b,0)$. En la Figure \ref{E35D}, vemos los caminos $\tb0\tb1\tb0\tb1\tb0\tb1$  y $\tb0\tb1\tb0\tb1\tb1\tb0\tb1\tb1$ respectivamente.\\

Alternativamente cada palabra puede ser codificada como un diagrama de Ferrers correspondiente este al conjunto de casillas a izquierda de la trayectoria. Es habitual, identificar a los diagramas Ferrers por el número de casillas en cada línea, lo que corresponde por lo tanto a las particiones:
 \begin{align}
 \l&=[\l_{a-1},\l_{a-2},\ldots,\l_{1}], &\text{  con  $ \l_{a-l}\leq\dint{\dfrac{bl}{a}}$ donde  $1\leq l\leq a-1$}. \label{Ferrers}
\end{align}

Por ejemplo en la Figura \ref{E35D}, los caminos son codificados  por las  secuencias $[2,1,0]$ and $[3,1,0]$, respectivamente. 
Los casos donde $(a,b)$ son coprimos, o $b=ak$ son de particular interes. Para el caso de $b=ak$, con $k\geq 1$, el cardinal esta dado por la formula de Fuss-Catalan (see \cite{KOS09}). 
\[
\Ct{(a,k)}=\dfrac{1}{ak+1}\cmb{ak+a}{a}.
\]
Si $a$ y $b$ son coprimos, obtemos la formula "general" de Catalan:
\[
\Ct{(a,b)}=\dfrac{1}{a+b}\cmb{a+b}{a}.
\]

En particular, cuando $a=p$ es  primo, $b$ y $p$ son coprimos, o $b$ es un multiplo de $p$. Por lo tanto el número correspondiente de caminos Dyck es:
\begin{align*}
|\D_{p,b}|=\begin{cases}
\frac{1}{p+b}\cmb{p+b}{p} &\mbox{ si } \tb{mcd}(p,b)=1, \\
\frac{1}{p+b+1}\cmb{p+b+1}{p} &\mbox{ si } b=kp.
\end{cases}
\end{align*} 

El problema boleta generalizado (ballot) se relaciona con el número de caminos de $ (0,0) $ hasta $ (a, b) $ inscriptos en un retículo que nunca superan la línea de $ y = kx $ (see \cite{Serrano03}):
\begin{align*}
\dfrac{b-ka+1}{b}\cmb{a+b}{a}&&\text{donde $k\geq 1$, y $b>ak\geq0$.}
\end{align*}
Ademas, el numero de caminos de longitud $2(k+1)n+1$ que comienzan en $(0,0)$ y evitan tocan o cruzar la  linea $y=kx$ (see \cite{Chap09}) esta dado por la formula:
\begin{align*} 
\cmb{2(k+1)n}{2n}-(k-1)\dsum{2n-1}{i=0}{\cmb{2(k+1)n}{i}},&&\text{donde $n\geq 1$ y $k\geq0$.}
\end{align*} 
En el caso mas general tenemos una formula debido a Bizley (see \cite{Biz54}) expresada como sigue.
Sea $m=da$, $n=db$ y $d=\tb{mcd}(m,n)$, entonces:
\begin{align*}
\Bz_{k}^{(a,b)}& :=\Biz{a}{b}{k} &\text{para $k\in \N$}, \\
\Bz_{\l}^{a,b}& :=\Bizu{a}{b}{\l} &\text{si $\l=(\l_1,\l_2,\dots,\l_l)$,}
\end{align*}

Es fácil demostrar que el número de caminos de Dyck en $m\times n$ $\D_{m,n}$ es:
\begin{align*}
|\D_{m,n}|&:= \dsum{\l \vdash d}{}{\frac{1}{z_{\l}}\Bz_{\l}^{(a,b)}}&\text{donde $n\geq 1$ y $k\geq0$.}
\end{align*}
\paragraph{\bf Palabra y camino de Christoffel.} Un camino de Christoffel entre dos puntos distintos $P=(0,k)$ y $P'=(0,l)$ en un retículo $a\times b$ es el camino mas cercano por debajo del segmento $PP'$ (see \cite{MR13}).  Pro ejemplo , el camino de Dyck de la Figure \ref{E35D}(b) es tambien Christoffel, y la palabra asociada se llama palabra de Christoffel. El camino de Christoffel de un retículo $a\times b$ es el camino de Christoffel asociado a la diagonal de dicho retículo.
 Como en el caso de los caminos de Dyck, cada camino de Christoffel en un retículo fijo $a\times b$ se identifica mediante un diagrama Ferrers de la forma $(\l_{a-1},\l_{a-2},\ldots,\l_{1})$ dada por la ecuación \eqref{Ferrers}.\\

\paragraph{\bf{Diagrama de Christoffel.}} Llamamos diagrama de Christoffel $\C_{a,b}$ al diagrama de Ferrers asociado al camino de Christoffel de un retículo $a \times b$. Éste es escrito como:
\begin{align}
 \l&=\left[\dint{\dfrac{(a-1)b}{a}},\dint{\dfrac{(a-2)b}{a}},\ldots,\dint{\dfrac{2b}{a}},\dint{\dfrac{b}{a}}\right].\label{Chris}
\end{align}\\

Para su uso posterior, definimos dos funciones asociadas al diagrama de Ferrers.
\begin{align}
Q_{a,b} &=\dfrac{(a-1)(b-1)+\tb{mcd}(a,b)-1}{2}.\label{Qab}
\end{align} 

Tambien, sea $ \De_ {a, b} (l) $ la diferencia entre las casillas de los diagramas de Ferrers diagrams asociados a los caminos de Christoffel de $ a\times b $ y $ a \times (b-1) $, respectivamente:

\begin{align}
\De_{a,b}(l):=\dint{\dfrac{bl}{a}}-\dint{\dfrac{(b-1)l}{a}},
\end{align}
donde $a<b\in\N$ y $1\leq l \leq a-1$.\\

\paragraph{\bf{Diagrama Isósceles.}} Llamamos diagrama isósceles $\I_{n}$ al diagrama de Ferrers asociado al camino de Christoffel de un retículo cuadrado de lado $n$. Dado un diagrama de Ferrers $\T_{\g}$, llamamos el \textit{diagrama isósceles máximo} al mas grande diagrama isósceles  incluido en $\T_{\g}$.\\

\paragraph{\bf{Conjunto de Ferrers.}} Sea $\T_{\g}$ el diagrama de Ferrers limitado por el camino $\g$. El conjunto de Ferrers de $\T_{\g}$ es el conjunto de todos loa caminos de  Dyck contenidos en $\T_{\g}$.\\

En las próximas secciones describimos métodos para calcular los caminos de ($a,b$)-Dyck cuando $a$ y $b$ no son necesariamente coprimos en términos de los números de Catalan. Por motivos de simplicidad, cada vez que mencionemos un diagrama isósceles estaremos pensando en un diagrama isósceles máximo.

\section{Método de comparación de diagramas de Ferrers}\label{FDCM}

Sea $\C_{a,b}$ le diagrama de Christoffel de $a\times b$. Para establecer los resultados principales necesitamos contar el numero de casillas en exceso  entre el diagrama del retículo  $a\times b$ y $a\times c$, con $c>b$. El metodo se desarrolla haciendo borrando las casillas en exceso entre ambos (ver \cite{BJE14}).
Usando las funciones $Q_{a,b}$ y $\De_{a,b}(l)$, obtenemos las siguientes reglas:
\begin{enumerate}[\tb{Regla} 1:]

\item \label{R1} Si $Q_{a,b}=1$ y $\De_{a,b}(i)=1$, hay una sola casilla en una esquina en $\C_{a,c}$ que no pertenece a $\C_{a,b}$. Sean $\T_{a_1,b_1}$ y $\T_{a'_1,b'_1}$ los diagramas de Ferrers  obtenidos a partir de borrar de $\C_{a,c}$ el rectángulo que contiene la casilla en exceso $\rg{ \a}$ (ver Figure \ref{DifDiag}(b)).

\begin{figure}[h]
\centering
\subcaptionbox{Comparación entre $\C_{a,b}$ y $\C_{a,c}$ }[0.5\linewidth]{
$
\ytableausetup{mathmode, boxsize=0.8em }
\begin{ytableau}
\none [^{_{\l_{n}}}\,\,\,\, ]&  \, &\,&\\
\none [^{_{\l_{n-1}}}\,\,\,\, ]&  \,& \, && \, &\\
\none [^{_{\l_{n-2}}} \,\,\,\,]&  \,& \, \,& \, & \,  & & \, & \\
\none [^{_{\vdots}} \,\,\,\,]& \none[^{_{\vdots}}]& \none[^{_{\vdots}}] \,& \none[^{_{\vdots}}] &
\none[^{_{\vdots}}]  &\none[^{_{\vdots}}] & \none[^{_{\vdots}}]& \none[^{_{\vdots}}] \\
\none [^{_{\l_{i}}}\,\,\,\, ]&  \,& \, \,& \, & \,  & \,& \, & & \, & &*(red)&\none[ \rg{ \a}] \\
\none [^{_{\vdots}} \,\,\,\,]& \none[^{_{\vdots}}]& \none[^{_{\vdots}}] \,& \none[^{_{\vdots}}] &
\none[^{_{\vdots}}]  &\none[^{_{\vdots}}] & \none[^{_{\vdots}}]& \none[^{_{\vdots}}] &\none[^{_{\vdots}}] & \none[^{_{\vdots}}]& \none[^{_{\vdots}}] \\
\none [^{_{\l_{2}}} \,\,\,\,]&   \,& \, \,& \, & \,  & \,& \, & \,&\, &\,& & \,&  \\
\none [^{_{\l_{1}}} \,\,\,\,]&   \,& \, \,& \, & \,  & \,& \, & \,&\, &\,& &\, &\, & & \\
\end{ytableau}
$} 
\subcaptionbox{$\T_{a_1,b_1}$ y $\T_{a'_1,b'_1}$}[0.3\linewidth]{
$
\ytableausetup{mathmode, boxsize=0.8em }
\begin{ytableau}
\none [^{_{\l_{n}}}\,\,\,\, ]&  \, &\,&\\
\none [^{_{\l_{n-1}}}\,\,\,\, ]&  \,& \, && \, &\\
\none [^{_{\l_{n-2}}} \,\,\,\,]&  \,& \, \,& \, & \,  & & \, & \\
\none [^{_{\vdots}} \,\,\,\,]& & &  &   &  & & &  & \\
\none [^{_{\l_{i}}}\,\,\,\, ]& \none[]& \none[] \,& \none[] & \none[]  & \none[]& \none[] & \none[]& \none[\cdots] &\none[] &*(red)  &\none[ \rg{ \a}] \\
\none [^{_{\vdots}} \,\,\,\,]& \none[]& \none[] \,& \none[] &\none[]  &\none[] & \none[]& \none[] &\none[] & \none[]& \none[^{_{\vdots}}] \\
\none [^{_{\l_{2}}} \,\,\,\,]&   \none[]& \none[] \,& \none[] & \none[]  & \none[]& \none[] & \none[]&\none[] &\none[]&\none[] & &  \\
\none [^{_{\l_{1}}} \,\,\,\,]&  \none[]& \none[] \,& \none[] & \none[]  & \none[]& \none[] & \none[]&\none[] &\none[]&\none[] && & & \\
\end{ytableau}
$
}
\caption{Regla 1.}\label{DifDiag}  
\end{figure}

Éstos diagramas de Ferrers no estan asociados en general a caminos de Christoffel. Sean \[\J_{a_1,b_1}\subseteq \D_{a_1,b_1}, \text{  y  } \J_{a'_1,b'_1}\subseteq \D_{a'_1,b'_1}\]  subconjuntos de los conjuntos de caminos de Dyck contenidos en los diagramas de Ferrers $\T_{a_1,b_1}$ y $\T_{a'_1,b'_1}$, respectivamente, entonces obtenemos:
\[|\D_{a,c}|-|\D_{a,b}|=-|\J_{a_1,b_1}|\cdot|\J_{a_1',b_1'}|,\]

Es claro ver que si la casilla $\a$ esta situada en la linea inferior del diagrama ($l=a-1$), la ecuación se reduce a:
\[|\D_{a,c}|-|\D_{a,b}|=-|\J_{a_1,b_1}|.\]
\item\label{R2} 
Cuando $Q_{a,b}=k$ y hay exactamente $k$ lineas differentes de una casilla de diferencia cada una necesitamos calcular cuantos caminos contienen dichas casillas (ver Figura \ref{DifDiag2}), para ello construimos una secuencia de conjuntos disjuntos como sigue. Sea $A_j $ el conjunto de caminos que no contiene las casillas $\rg{\a_i}$ para $ i> j $, donde $1\leq j\leq k$.  Tambien podemos construir en sentido inverso, sea $B_j$ los conjuntos de caminos que no contienen las casillas $ \a_i$ para $i <j $, donde $ 1 \leq j \leq k $.

\begin{figure}[h]
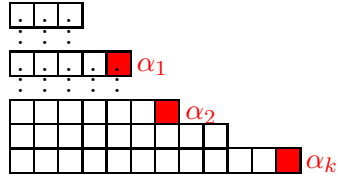

$$
\ytableausetup{mathmode, boxsize=0.8em}
\begin{ytableau}
\none [  ]&  \, &&\\
\none [ ]&\none [^{_{\vdots}}] &\none [^{_{\vdots}}]&\none [^{_{\vdots}}]\\
\none []&  \,& && \ &*(red) & \none[\rg{\ \ \a_1}]\\
\none [ ]&\none [^{_{\vdots}}] &\none [^{_{\vdots}}]&\none [^{_{\vdots}}] &\none [^{_{\vdots}}]&\none [^{_{\vdots}}]\\
\none []&  \,&  \,&  &   & &  & *(red) & \none[\ \ \rg{\a_2}] \\
\none [ ]&  \,& \,&  &   & &  & & &  \\
\none []&   \,&   &   & && & && & & &*(red)& \none[\ \ \rg{\a_k}] \\
\end{ytableau}
$$ 
\caption{Mas de una casilla.}\label{DifDiag2}
\end{figure}

Es facil ver que esta estrategia los provee de dos secuencias de conjuntos disjountos dos a dos tales que preservan la unión total, usando la \tb{Regla \ref{R1}} para cada $A_j$ o $B_j$ obtenemos:
\[|\D_{a,c}|-|\D_{a,b}|=-\dsum{j=1}{k}{(|\J_{a_j,b_j}|\cdot|\J_{a_j',b_j'}|)},\]

donde $\J_{a_j,b_j}\subseteq \D_{a_j,b_j}$, and $\J_{a'_j,b'_j}\subseteq \D_{a'_j,b'_j}$.
\end{enumerate}

\section{Método de descomposición de Diagramas}
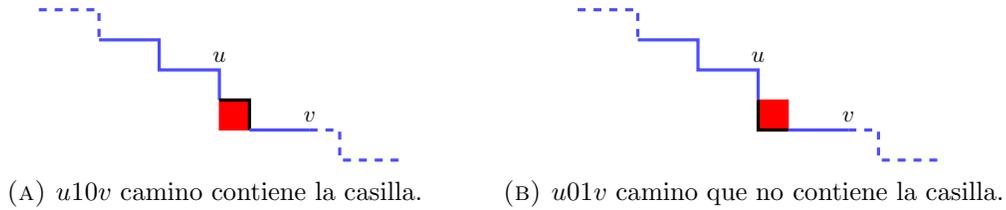
\begin{figure}[h]
\centering
\subcaptionbox{$u10v$ camino contiene la casilla.}[0.4\linewidth]{
\begin{tikzpicture}[scale=.4, every node/.style={scale=0.8}] 
\draw[fill=red,draw=red] ( 2,-1 ) -- ( 2,-2 ) -- (3,-2) -- (3,-1)-- cycle; 
\draw[color=blue!70,dashed,line width=1.2pt](-4,2)--(-2,2)--(-2,1);
\draw[color=blue!70,dashed,line width=1.2pt](5,-2)--(6,-2)--(6,-3)--(8,-3);
\draw[color=blue!70,line width=1.2pt](-2,1)--(0,1)--(0,0) -- ( 2, 0 ) -- ( 2,-1 );
\draw[color=black,line width=1.2pt] ( 2,-1 )--(3,-1)--(3,-2);
\draw[color=blue!70,line width=1.2pt](3,-2)--(5,-2); 
\draw [ black ]  ( 2,0 )  node[above] {$u$ };
\draw [ black ]  ( 5,-2 )  node[above] {$v$ };
\end{tikzpicture}
}
\subcaptionbox{$u01v$ camino que no contiene la casilla.}[0.4\linewidth]{
\begin{tikzpicture}[scale=.4, every node/.style={scale=0.8}] 
\draw[fill=red,draw=red] ( 2,-1 ) -- ( 2,-2 ) -- (3,-2) -- (3,-1)-- cycle; 
\draw[color=blue!70,dashed,line width=1.2pt](-4,2)--(-2,2)--(-2,1);
\draw[color=blue!70,dashed,line width=1.2pt](5,-2)--(6,-2)--(6,-3)--(8,-3);
\draw[color=blue!70,line width=1.2pt](-2,1)--(0,1)--( 0, 0 ) -- ( 2, 0 ) -- ( 2,-1 );
\draw[color=black,line width=1.2pt]( 2,-1 )--(2,-2)--(3,-2); 
\draw[color=blue!70,line width=1.2pt](3,-2)--(5,-2); 
\draw [ black ]  ( 2,0 )  node[above] {$u$ };
\draw [ black ]  ( 5,-2 )  node[above] {$v$ };
\end{tikzpicture}
}
\caption{Separación de diagramas.}\label{suma}
\end{figure}

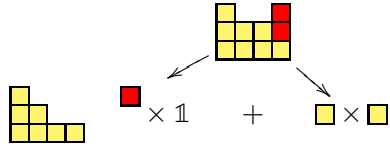
\begin{figure}[h] 
\begin{align*}
\xymatrixrowsep{0.05in}
\xymatrixcolsep{0.05in}
\xymatrix{
&\ytableausetup {mathmode, boxsize=0.6em,centertableaux} 
\ydiagram[*(red)]
  {1+0,2+1,3+1}
*[*(yellow!70)]{1,3,4}
\ar@{->}[ld]^{} \ar@{->}[rd]^{}  \\
\ytableausetup {mathmode, boxsize= 0.6em,centertableaux} 
\ydiagram[*(red)]
  {1+0,2+0,3+1}
*[*(yellow!70)]{1,2,4}
\times \U& +  &\ytableausetup {mathmode, boxsize= 0.6em,centertableaux} \ydiagram [ *(yellow!70)]{ 1 }  \times \ydiagram [ *(yellow!70)]{ 1 }& \\
}
\end{align*}
\caption{Primera iteración en la descomposición.}\label{pas1}
\end{figure}

Usando el método de comparación (\tb{Regla 1}), podemos iterar el proceso de borrado de casillas en exceso entre un diagrama cualquiera  $\T_{\g}$ y su respectivo diagrama isósceles $\I_{n}$. Comenzando por la casilla superior derecha (siempre de una esquina) como se muestra en la  Figura \ref{pas1}. Dicha descomposición esta dada en sumas y productos. La operación de suma $+$ esta dad por la unión de conjuntos de Ferrers disjuntos. Consideremos una casilla (roja) en el borde del diagrama. Una camino contenido puede ser escrito como se ve en la figura \ref{suma}.
El producto de diagramas $\T\times\T'$ es un diagrama conteniendo todas las posible concatenaciones  entre un camino de Dick de $\T$ y uno de $\T'$. Ambas operaciones relacionan los diagramas con el cardinal del cardinal de los conjuntos de Ferrers. Vale aclarar que los diagramas isosceles de lado $n$ contienen $\Ct{n}$ numero de caminos, destacando que el diagrama vacío $\U$ contiene un solo camino.\\
La relación $\H$ entre los diagramas y los números naturales $\N$ es definida de forma que si consideramos a los diagramas isosceles obtenemos: 
\begin{align*}
\H(\U)&:=1,&\\
\H(\I_n)&:=\Ct{n},&\\
\H(\I_n+\I_m)&:=\H(\I_n)+\H(\I_m),&\\
\H(\I_n\times\I_m)&:=\H(\I_n)\cdot\H(\I_m).&
\end{align*}
Debido a la definición de las operaciones podemos extender estos resultados a cualquier diagrama, permitiéndonos calcular el cardinal de cualquier conjunto de Ferrers. 

Por ejemplo, el diagrama de Christofel  $\C_{4,6}$ es $[4,3,1]$ $\ytableausetup {mathmode, boxsize= 0.5 em,centertableaux} \ydiagram [ ]{ 1 , 3 , 4 }$ y su correspondiente isósceles es $[3,2,1]$ $\ytableausetup {mathmode, boxsize= 0.5 em,centertableaux} \ydiagram [*(yellow!70) ]{ 1 , 2 , 3 }$. cuando borramos la primera casilla del diagrama original se divide en dos pares de diagramas asociados con las operaciones que simplifican el calculo de los caminos (see Figure \ref{pas1}).\\

En la Figura \ref{pas1}, $\U$ es el diagrama vacío. Podemos repetir este método hasta que todos los diagramas sean isósceles (la operación $\times$ es distributiva respecto a la operación $+$).
\begin{figure}[h] 
\begin{align*}
\xymatrixrowsep{0.03in}
\xymatrixcolsep{0.03in}
\xymatrix{
&&\ytableausetup {mathmode, boxsize= 0.6em,centertableaux}
\ydiagram[*(red)]
  {1+0,2+1,3+1}
*[*(yellow!70)]{1,3,4}
\ar@{->}[ld]^{} \ar@{->}[rd]^{}  \\
&\ytableausetup {mathmode, boxsize= 0.6em,centertableaux}
\ydiagram[*(red)]
  {1+0,2+0,3+1}
*[*(yellow!70)]{1,2,4}\times \U\ar@{->}[ld]^{}\ar@{->}[rd]^{}& +  &\ytableausetup {mathmode, boxsize= 0.6em ,centertableaux} \ydiagram [ *(yellow!70)]{ 1 }  \times \ydiagram [ *(yellow!70)]{ 1 } \ar@{->}[rd]^{}& \\
\ytableausetup {mathmode, boxsize= 0.6em,centertableaux} \ydiagram [ *(yellow!70)]{ 1 , 2 , 3 } \times \U\times \U&+& \ydiagram [ *(yellow!70)]{ 1 , 2  }\times \U\times \U &+&\ydiagram [ *(yellow!70)]{ 1 }  \times \ydiagram [ *(yellow!70)]{ 1 }&&\\
}
\end{align*}
\caption{Descomposición completa del diagrama.}\label{dec2}
\end{figure}
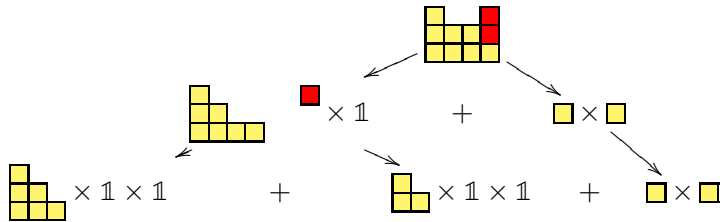

\section{Método alternativo de descomposición}
\begin{figure}[h]
\begin{align*}
\xymatrixrowsep{0.05in}
\xymatrixcolsep{0.05in}
\xymatrix{
&\ytableausetup {mathmode, boxsize=0.6em,centertableaux} 
\ydiagram[*(green)]
  {1+0,2+0,2+1}
*[*(yellow!70)]{1,2,3}
\ar@{->}[ld]^{} \ar@{->}[rd]^{}  \\
\ytableausetup {mathmode, boxsize= 0.6em,centertableaux} 
\ydiagram[*(red)]
  {1+0,2+0,2+0}
*[*(yellow!70)]{1,2,2}
\times \U& +  &\ytableausetup {mathmode, boxsize= 0.6em,centertableaux} \ydiagram[*(red)]
  {1+0,2+0}
*[*(yellow!70)]{1,2} \times\U& \\
}
\end{align*}
\caption{Descomposición alternativa.}\label{decal}
\end{figure}
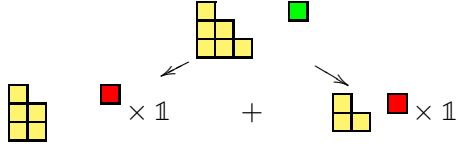

Podemos utilizar el método de descomposición de diagramas ya descripto, pero esta vez comparando las casillas entre el diagrama $\T_{\g}$ y el diagrama isosceles $\I_n$ mas pequeño que lo contenga. De esta manera obtendremos un desarrollo del diagrama $\I_n$ en cual aparezca el diagrama que deseamos calcular. Por ejemplo, para el  [2,2,1] obtenemos la siguente descomposición (ver Figura \ref{decal}).

Entonces,
\begin{align*}
\H\left(\ytableausetup {mathmode, boxsize= 0.6 em,centertableaux} \ydiagram [
*(yellow!70) ]{ 1,2 , 3 } 
\right)&=\H\left(\ytableausetup {mathmode, boxsize= 0.6 em,centertableaux} \ydiagram [
*(yellow!70) ]{ 1,2, 2 }\times\U + \ydiagram [
*(yellow!70) ]{ 1, 2 }\times\U\right)&\\
&=\H\left(\ytableausetup {mathmode, boxsize= 0.6 em,centertableaux} \ydiagram [
*(yellow!70) ]{ 1,2, 2 }\right) + \H\left(\ydiagram [
*(yellow!70) ]{ 1, 2 }\right)&\\
\Ct{4}&=\H\left(\ytableausetup {mathmode, boxsize= 0.6 em,centertableaux} \ydiagram [
*(yellow!70) ]{ 1,2, 2 }\right)+\Ct{3}&\\
&\Rightarrow&\\
\H\left(\ytableausetup {mathmode, boxsize= 0.6 em,centertableaux} \ydiagram [
*(yellow!70) ]{ 1,2, 2 }\right)&=\Ct{4}-\Ct{3}=9&\\
\end{align*}

\section{Teorema, lema y corolario}
\begin{teorema}\label{t1}
Dado un diagrama de Ferrers $\T_{\g}$ siempre puede ser divido en dos pares de diagramas uno resultado de haber borrado una casilla de su borde $\tau$ multiplicado por el diagrama vacio y el otro resultado de haber borrado el rectángulo inferior izquierdo que contiene la casilla $\tau$, de forma que no alteren el cardinal del conjunto de Ferrers de $\T_{\g}$.
\end{teorema}
\begin{prueva}
\

Utilizando el Lema \ref{l1} queda clara la forma de separación.
\end{prueva}

\begin{lema}\label{l1} Sea $\T_{\g}$ un diagrama de Ferrers limitado por un camino de Dyck $\g$ de $\D_{a,b}$. Sea $\tau$ una casilla del retículo situada en el borde de nuestro camino $\g$ de coordenadas $(c,d)$, consideramos los diagramas $\Tx_{\g}$ y $\Ty_{\g}$, el primero resultado de haber borrado la casilla $\tau$ al diagrama $\T_{\g}$ y el segundo resultado de haber borrado el rectngulo inferior izquierdo que contiene a la casilla $\tau$ de $\T_{\g}$ (ver figura \ref{DifDiag}). 
\begin{enumerate}[a$)$]
\item los conjuntos $\Tx_{\g}$ y $\Ty_{\g}$ son disjuntos, 
\item para todo $\de \in\Ty_{\g}$ existen  $\a$ y $\b$ fijos tales que $\de=\e_1+\e_2$ donde $\e_1\in \T_{\a}$ y  $\e_2\in \T_{\b}$,
\item sean $|\Ty_{\g}|=|\T_{\a}|\cdot|\T_{\b}|$.
\end{enumerate}
\end{lema}

\begin{prueva}
\
\begin{enumerate}[a$)$]
\item sea $\de_1\in\Tx_{\g}$ y $\de_2\in\Ty_{\g}$, sean $w_{\de_1}$ y $w_{\de_2}$ las palabras asociadas respectivamente. En el caso general son distintos estudiaremos en el caso donde podrían no serlo. Como  $\de_2\in\Ty_{\g}$, $w_{\de_2}=u10v$ tel que $|u|_0=d$, $|u|_1=c-1$, $|v|_0=(a-c)-1$ y $|v|_1=(b-d)$, si $w_{\de_1}=usv$ tal que $|s|_0=|s|_1=1$, entonces $s\neq 10$ ya que por definición de $\Tx_{\g}$ $\de_1$ no contiene la casilla $\tau$, luego ambos conjuntos son disjuntos.

\item como $\g \in\Ty_{\g}$, $w_{\g}=u10v$ sean $\a$ y $\b$ los caminos tales que $w_{\a}=u1$ y $w_{\b}=0v$, donde $u$ y $v$ tienen las mismas condiciones que en el item a), obtenemos que para todo  $\de\in\Ty_{\g}$, $\de=\e_1+\e_2$ donde $\e_1\in \T_{\a}$ y  $\e_2\in \T_{\b}$,
\item Como $\g=\a+\b$ y por cada camino que pertenece a $\de\in\T_{\a}$ podemos concatenarlo con cualquier camino $\e\in\T_{\b}$ para obtener un camino en $\T_{\g}$, luego $|\Ty_{\g}|=|\T_{\a}|\cdot|\T_{\b}|$.
\end{enumerate}

\end{prueva}


\begin{corolario}
Para todo diagrama de Ferrers $\T_{\g}$, utilizando la separación según el Teorema \ref{t1} y la relación $\H$ obtenemos la cantidad de caminos contenidos $|\T_{\g}|$ 
\end{corolario}
\begin{prueva}
\

Utilizando el Teorema \ref{t1} en forma recursiva hasta obtener solo diagramas isosceles y por el Lema \ref{l1} queda claro que cada separación nos da conjuntos disjuntos con respecto a los caminos contenidos, luego aplicando $\H$ obtenemos el numero total de caminos.
\end{prueva}



\section{Ejemplos}
En esta sección desarrollamos varios ejemplos para mostrar la aplicación del método.
\subsection{Desarrollo de [1,1]}

\begin{align*}
\xymatrixrowsep{0.05in}
\xymatrixcolsep{0.05in}
\xymatrix{
&\ytableausetup {mathmode, boxsize=0.6em,centertableaux} 
\ydiagram[*(red)]{0+1,  1+0} *[*(yellow!70)]{1,1}
\ar@{->}[ld]^{} \ar@{->}[rd]^{}  \\
\ytableausetup {mathmode, boxsize= 0.6em,centertableaux} 
\ydiagram[*(red)]{1+0} *[*(yellow!70)]{1}
\times \U& +  &\ytableausetup {mathmode, boxsize= 0.6em,centertableaux} \U\times\U& \\
}
\end{align*}

\begin{align*}
\H\left(\ytableausetup {mathmode, boxsize= 0.6em,centertableaux} \ydiagram [*(yellow!70)]{1,1} \right)&=\H\left(\ytableausetup {mathmode, boxsize= 0.6em,centertableaux} \ydiagram [
*(yellow!70) ]{ 1 } \times \U + \U \times \U\right)&\\
&=\H\left(\ytableausetup {mathmode, boxsize= 0.6 em,centertableaux} \ydiagram [
*(yellow!70) ]{ 1 }\times\U\right) + \H\left(\U\times\U\right) &\\
&=\H\left(\ytableausetup {mathmode, boxsize= 0.6 em,centertableaux} \ydiagram [
*(yellow!70) ]{ 1 }\right)\cdot\H\left(\U\right) + \H\left(\U\right)\cdot\H\left(\U\right) &\\
&=\H\left(\I_2\right)\cdot\H\left(\U\right) + \H\left(\U\right)\cdot\H\left(\U\right) &\\
&=\Ct{2}+1&\\
&=2+1=3&
\end{align*}

\subsection{Desarrollo de [2,2]}

\begin{align*}
\xymatrixrowsep{0.05in}
\xymatrixcolsep{0.05in}
\xymatrix{
&\ytableausetup {mathmode, boxsize=0.6em,centertableaux} 
\ydiagram[*(red)]
  {1+1,2+0}
*[*(yellow!70)]{2,2}
\ar@{->}[ld]^{} \ar@{->}[rd]^{}  \\
\ytableausetup {mathmode, boxsize= 0.6em,centertableaux} 
\ydiagram[*(red)]
  {1+0,2+0}
*[*(yellow!70)]{1,2}
\times \U& +  &\ytableausetup {mathmode, boxsize= 0.6em,centertableaux} \U\times\U& \\
}
\end{align*}

\begin{align*}
\H\left(\ytableausetup {mathmode, boxsize= 0.6 em,centertableaux} \ydiagram [
*(yellow!70) ]{ 2 , 2 } 
\right)&=\H\left(\ytableausetup {mathmode, boxsize= 0.6 em,centertableaux} \ydiagram [
*(yellow!70) ]{ 1 , 2 }\times\U + \U\times\U\right) &\\
&=\H\left(\ytableausetup {mathmode, boxsize= 0.6 em,centertableaux} \ydiagram [
*(yellow!70) ]{ 1 , 2 }\times\U\right) + \H\left(\U\times\U\right) &\\
&=\H\left(\ytableausetup {mathmode, boxsize= 0.6 em,centertableaux} \ydiagram [
*(yellow!70) ]{ 1 , 2 }\right)\cdot\H\left(\U\right) + \H\left(\U\right)\cdot\H\left(\U\right) &\\
&=\H\left(\I_3\right)\cdot\H\left(\U\right) + \H\left(\U\right)\cdot\H\left(\U\right) &\\
&=\Ct{3}+1&\\
&=5+1=6&
\end{align*}

\subsection{Desarrollo de [2,2,1]}\label{e221}

\begin{align*}
\xymatrixrowsep{0.05in}
\xymatrixcolsep{0.05in}
\xymatrix{
&\ytableausetup {mathmode, boxsize=0.6em,centertableaux} 
\ydiagram[*(red)]
  {0+1,1+1,2+0}
*[*(yellow!70)]{1,2,2}
\ar@{->}[ld]^{} \ar@{->}[rd]^{}  \\
\ytableausetup {mathmode, boxsize= 0.6em,centertableaux} 
\ydiagram[*(red)]
  {1+1,2+0}
*[*(yellow!70)]{2,2}
\times \U& +  &\ytableausetup {mathmode, boxsize= 0.6em,centertableaux} \U\times\ydiagram[*(red)]
  {0+1,1+0}
*[*(yellow!70)]{1,1}& \\
}
\end{align*}

\begin{align*}
\H\left(\ytableausetup {mathmode, boxsize= 0.6 em,centertableaux} \ydiagram [
*(yellow!70) ]{ 1,2 , 2 } 
\right)&=\H\left(\ytableausetup {mathmode, boxsize= 0.6 em,centertableaux} \ydiagram [
*(yellow!70) ]{ 2, 2 }\times\U + \U\times\ydiagram [
*(yellow!70) ]{ 1, 1 }\right) &\\
&=\H\left(\ytableausetup {mathmode, boxsize= 0.6 em,centertableaux} \ydiagram [
*(yellow!70) ]{ 2 , 2 }\times\U\right) + \H\left(\U\times\ydiagram [*(yellow!70) ]{ 1, 1 }\right)&\\
&=\H\left(\ytableausetup {mathmode, boxsize= 0.6 em,centertableaux} \ydiagram [
*(yellow!70) ]{ 2 , 2 }\right)\cdot\H\left(\U\right) + \H\left(\U\right)\cdot\H\left(\ydiagram [*(yellow!70) ]{ 1, 1 }\right) &\\
&=\Ct{3}+1+\Ct{2}+1&\\
&=9&
\end{align*}

\subsection{Desarrollo de [3,1,1]}

\begin{align*}
\xymatrixrowsep{0.05in}
\xymatrixcolsep{0.05in}
\xymatrix{
&&\ytableausetup {mathmode, boxsize=0.6em,centertableaux} 
\ydiagram[*(red)]
  {0+1,1+0,2+1}
*[*(yellow!70)]{1,1,3}
\ar@{->}[ld]^{} \ar@{->}[rd]^{}  \\
&\ytableausetup {mathmode, boxsize= 0.6em,centertableaux} 
\ydiagram[*(red)]
  {1+0,2+1}
*[*(yellow!70)]{1,3}
\times \U \ar@{->}[ld]^{}\ar@{->}[rd]^{} & +  &\ytableausetup {mathmode, boxsize= 0.6em,centertableaux} \U\times\ydiagram[*(red)]
  {1+1}
*[*(yellow!70)]{2}& \\
\ytableausetup {mathmode, boxsize= 0.6em,centertableaux} 
\ydiagram[*(red)]
  {1+0,2+0}
*[*(yellow!70)]{1,2}
\times \U \times \U &+&\ytableausetup {mathmode, boxsize= 0.6em,centertableaux} 
\ydiagram[*(red)]
  {1+0}
*[*(yellow!70)]{1}
\times \U \times \U\\
}
\end{align*}

\begin{align*}
\H\left(\ytableausetup {mathmode, boxsize= 0.6 em,centertableaux} \ydiagram [
*(yellow!70) ]{ 1,1 , 3 } 
\right)&=\H\left(\ytableausetup {mathmode, boxsize= 0.6 em,centertableaux} \ydiagram [
*(yellow!70) ]{ 1, 2 }\times\U \times\U +\ydiagram [
*(yellow!70) ]{ 1 }\times\U \times\U+ \U\times\ydiagram [
*(yellow!70) ]{ 2 }\right) &\\
&=\H\left(\ytableausetup {mathmode, boxsize= 0.6 em,centertableaux} \ydiagram [
*(yellow!70) ]{ 1, 2 }\right) +\H\left(\ydiagram [
*(yellow!70) ]{ 1 }\right)+ \H\left(\ydiagram [
*(yellow!70) ]{ 2 }\right) &\\
&=\H\left(\ytableausetup {mathmode, boxsize= 0.6 em,centertableaux} \ydiagram [
*(yellow!70) ]{ 1 , 2 }\right) +\H\left(\ydiagram [*(yellow!70) ]{ 1}\right)+ \H\left(\ydiagram [*(yellow!70) ]{ 1,1 }\right)&\\
&=\Ct{3}+\Ct{2}+\Ct{2}+1&\\
&=\Ct{3}+2\Ct{2}+1&\\
&=5+4+1=10&
\end{align*}

\subsection{Desarrollo alternativo de [3,1,1]}\label{311a}

\begin{align*}
\xymatrixrowsep{0.05in}
\xymatrixcolsep{0.05in}
\xymatrix{
&&\ytableausetup {mathmode, boxsize=0.6em,centertableaux} 
\ydiagram[*(green)]
  {1+0,1+1,3+0}
*[*(yellow!70)]{1,2,3}
\ar@{->}[ld]^{} \ar@{->}[rd]^{}  \\
&\ytableausetup {mathmode, boxsize= 0.6em,centertableaux} 
\ydiagram[*(red)]
  {1+0,1+0,3+0}
*[*(yellow!70)]{1,1,3}
\times \U & +  &\ytableausetup {mathmode, boxsize= 0.6em,centertableaux}\ydiagram[*(yellow!70)]{1}\times\ydiagram[*(yellow!70)]{1}& 
}
\end{align*}

\begin{align*}
\H\left(\ytableausetup {mathmode, boxsize= 0.6 em,centertableaux} \ydiagram [
*(yellow!70) ]{ 1,2 , 3 } 
\right)&=\H\left(\ytableausetup {mathmode, boxsize= 0.6 em,centertableaux} \ydiagram [
*(yellow!70) ]{ 1, 1,3 }\times\U +\ydiagram [
*(yellow!70) ]{ 1 }\times\ydiagram [
*(yellow!70) ]{ 1 }\right) &\\
&=\H\left(\ytableausetup {mathmode, boxsize= 0.6 em,centertableaux} \ydiagram [
*(yellow!70) ]{ 1, 1,3 }\right) +\H\left(\ydiagram [
*(yellow!70) ]{ 1 }\right)\cdot\H\left(\ydiagram [
*(yellow!70) ]{ 1 }\right) &\\
\Ct{4}&=\H\left(\ytableausetup {mathmode, boxsize= 0.6 em,centertableaux} \ydiagram [
*(yellow!70) ]{ 1, 1,3 }\right) +\Ct{2}^2 &\\
&\Rightarrow&\\
\H\left(\ytableausetup {mathmode, boxsize= 0.6 em,centertableaux} \ydiagram [
*(yellow!70) ]{ 1, 1,3 }\right) &=\Ct{4}-\Ct{2}^2=10&\\
\end{align*}

\subsection{Desarrollo de [3,2,2]}

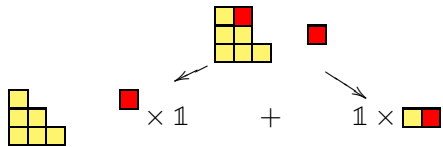
\begin{align*}
\xymatrixrowsep{0.05in}
\xymatrixcolsep{0.05in}
\xymatrix{
&\ytableausetup {mathmode, boxsize=0.6em,centertableaux} 
\ydiagram[*(red)]
  {1+1,2+0,3+0}
*[*(yellow!70)]{2,2,3}
\ar@{->}[ld]^{} \ar@{->}[rd]^{}  \\
\ytableausetup {mathmode, boxsize= 0.6em,centertableaux} 
\ydiagram[*(red)]
  {1+0,2+0,3+0}
*[*(yellow!70)]{1,2,3}
\times \U& +  &\ytableausetup {mathmode, boxsize= 0.6em,centertableaux} \U\times\ydiagram[*(red)]
  {1+0}
*[*(yellow!70)]{1}& \\
}
\end{align*}

\begin{align*}
\H\left(\ytableausetup {mathmode, boxsize= 0.6 em,centertableaux} \ydiagram [
*(yellow!70) ]{ 2,2 , 3 } 
\right)&=\H\left(\ytableausetup {mathmode, boxsize= 0.6 em,centertableaux} \ydiagram [
*(yellow!70) ]{ 1,2, 3 }\times\U + \U\times\ydiagram [
*(yellow!70) ]{ 1}\right) &\\
&=\H\left(\ytableausetup {mathmode, boxsize= 0.6 em,centertableaux} \ydiagram [
*(yellow!70) ]{ 1,2, 3 }\right) + \H\left(\ydiagram [
*(yellow!70) ]{ 1}\right) &\\
&=\Ct{4}+\Ct{2}&\\
&=14+2=16&
\end{align*}
\subsection{Desarrollo de [4,3,3,1]}
\begin{align*}
\ytableausetup {mathmode, boxsize= 0.6em,centertableaux} \ydiagram[*(red)]
  {1+0,2+1}*[
*(yellow!70) ]{ 1,3 , 3 , 4 }
= &\ytableausetup {mathmode, boxsize= 0.6em,centertableaux} \ydiagram [
*(yellow!70) ]{ 1 , 2 , 3 , 4 } \times \U + \ytableausetup {mathmode,
boxsize= 0.6em,centertableaux} \ydiagram [ *(yellow!70) ]{ 1 } \times
\ytableausetup {mathmode, boxsize= 0.6em,centertableaux} \ydiagram [
*(yellow!70) ]{ 1 }
\end{align*}

\begin{align*}
\H\left(\ytableausetup {mathmode, boxsize= 0.6em,centertableaux} \ydiagram [
*(yellow!70) ]{ 1, 3 , 3 , 4 }
\right)&=\Ct{5}+\Ct{2}\cdot\Ct{2}&\\
&=42+4=46&
\end{align*}
\subsection{Desarrollo de [5,4,2,1]}

\begin{align*}
\ytableausetup {mathmode, boxsize= 0.6em,centertableaux} \ydiagram[*(red)]
  {1+0,2+0,3+1,4+1}*[
*(yellow!70) ]{ 1,2 , 4 , 5 }
= & \ytableausetup {mathmode, boxsize= 0.6em,centertableaux} \ydiagram [
*(yellow!70) ]{ 1 , 2 , 3 , 5 } \times \U + \ytableausetup {mathmode,
boxsize= 0.6em,centertableaux} \ydiagram [ *(yellow!70) ]{ 1 } \times
\ytableausetup {mathmode, boxsize= 0.6em,centertableaux} \ydiagram [
*(yellow!70) ]{ 1 , 2 } &\\
=& \U \times \ytableausetup {mathmode, boxsize= 0.6em,centertableaux}
\ydiagram [ *(yellow!70) ]{ 1 , 2 , 3 , 4 } \times \U + \U \times \U \times
\ytableausetup {mathmode, boxsize= 0.6em,centertableaux} \ydiagram [
*(yellow!70) ]{ 1 , 2 , 3 } + \ytableausetup {mathmode, boxsize=0.6em,centertableaux} \ydiagram [ *(yellow!70) ]{ 1 } \times \ytableausetup
{mathmode, boxsize= 0.6em,centertableaux} \ydiagram [ *(yellow!70) ]{ 1 ,
2 } &\\
=&\ytableausetup {mathmode, boxsize= 0.6em,centertableaux} \ydiagram [
*(yellow!70) ]{ 1 , 2 , 3 , 4 } + \ytableausetup {mathmode, boxsize= 0.6em,centertableaux} \ydiagram [ *(yellow!70) ]{ 1 , 2 , 3 } +
\ytableausetup {mathmode, boxsize= 0.6em,centertableaux} \ydiagram [
*(yellow!70) ]{ 1 } \times \ytableausetup {mathmode, boxsize= 0.6em,centertableaux} \ydiagram [ *(yellow!70) ]{ 1 , 2 }
\end{align*}
\begin{align*}
\H\left(\ytableausetup {mathmode, boxsize= 0.6em,centertableaux} \ydiagram [
*(yellow!70) ]{ 1, 2 , 4 , 5 }
\right)&=\Ct{5}+\Ct{4}+\Ct{3}\cdot\Ct{2}&\\
&=42+14+10=66&
\end{align*}
\subsection{Développe de $[5,5,2,2]$}\label{it8b}
\

\

Para desarrollar este diagrama hacen falta 8 iteraciones después de las cuales el diagrama queda reducido y transformado a un polinomio en números de Catalan como sigue:
\begin{align*}
\ytableausetup {mathmode, boxsize= 0.6em,centertableaux}\ydiagram[*(red)]
  {1+1,2+0,3+2,4+1}*[
*(yellow!70) ]{ 2,2 , 5 , 5 }
= & \ytableausetup {mathmode, boxsize= 0.6em,centertableaux} \ydiagram [
*(yellow!70) ]{ 1 , 2 , 5 , 5 } \times \U + \ytableausetup {mathmode,
boxsize= 0.6em,centertableaux} \ydiagram [ *(yellow!70) ]{ 3 , 3 } \times
\U &\\
=&\U \times \ytableausetup {mathmode, boxsize= 0.6em,centertableaux}
\ydiagram [ *(yellow!70) ]{ 1 , 2 , 4 , 5 } \times \U + \U \times \U \times
\ytableausetup {mathmode, boxsize= 0.6em,centertableaux} \ydiagram [
*(yellow!70) ]{ 1 , 2 } + \U \times \ytableausetup {mathmode, boxsize= 0.6em,centertableaux} \ydiagram [ *(yellow!70) ]{ 2 , 3 } \times \U + \U
\times \U \times \U &\\
=&\U \times \U \times \ytableausetup {mathmode, boxsize= 0.6em,centertableaux} \ydiagram [ *(yellow!70) ]{ 1 , 2 , 3 , 5 } \times \U
+ \U \times \U \times \ytableausetup {mathmode, boxsize= 0.6em,centertableaux} \ydiagram [ *(yellow!70) ]{ 1 } \times \ytableausetup
{mathmode, boxsize= 0.6em,centertableaux} \ydiagram [ *(yellow!70) ]{ 1 ,
2 } + \U \times \U \times \ytableausetup {mathmode, boxsize= 0.6em,centertableaux} \ydiagram [ *(yellow!70) ]{ 1 , 2 } + &\\
&+\U \times \U\times \ytableausetup {mathmode, boxsize= 0.6em,centertableaux} \ydiagram
[ *(yellow!70) ]{ 1 , 3 } \times \U + \U \times \U \times \ytableausetup
{mathmode, boxsize= 0.6em,centertableaux} \ydiagram [ *(yellow!70) ]{ 1 }
\times \U + \U \times \U \times \U&\\
=&\U \times \U \times \U \times \ytableausetup {mathmode, boxsize= 0.6em,centertableaux} \ydiagram [ *(yellow!70) ]{ 1 , 2 , 3 , 4 } \times \U
+ \U \times \U \times \U \times \U \times \ytableausetup {mathmode, boxsize=
0.6em,centertableaux} \ydiagram [ *(yellow!70) ]{ 1 , 2 , 3 } +&\\
&+ \U \times\U \times \ytableausetup {mathmode, boxsize= 0.6em,centertableaux}
\ydiagram [ *(yellow!70) ]{ 1 } \times \ytableausetup {mathmode,
boxsize= 0.6em,centertableaux} \ydiagram [ *(yellow!70) ]{ 1 , 2 } + \U
\times \U \times \ytableausetup {mathmode, boxsize= 0.6em,centertableaux}
\ydiagram [ *(yellow!70) ]{ 1 , 2 } + \U \times \U \times \U \times
\ytableausetup {mathmode, boxsize= 0.6em,centertableaux} \ydiagram [
*(yellow!70) ]{ 1 , 2 } \times \U +&\\
&+ \U\times \U \times \U \times \U \times
\ytableausetup {mathmode, boxsize=0.6em,centertableaux} \ydiagram [
*(yellow!70) ]{ 1 } + \U \times \U \times \ytableausetup {mathmode,
boxsize=0.6em,centertableaux} \ydiagram [ *(yellow!70) ]{ 1 } \times \U +
\U \times \U \times \U&\\
=&\ytableausetup {mathmode, boxsize=0.6em,centertableaux} \ydiagram [
*(yellow!70) ]{ 1 , 2 , 3 , 4 } + \ytableausetup {mathmode, boxsize=0.6em,centertableaux} \ydiagram [ *(yellow!70) ]{ 1 , 2 , 3 } +
\ytableausetup {mathmode, boxsize=0.6em,centertableaux} \ydiagram [
*(yellow!70) ]{ 1 } \times \ytableausetup {mathmode, boxsize=0.6em,centertableaux} \ydiagram [ *(yellow!70) ]{ 1 , 2 } + \ytableausetup
{mathmode, boxsize=0.6em,centertableaux} \ydiagram [ *(yellow!70) ]{ 1 ,
2 } + \ytableausetup {mathmode, boxsize=0.6em,centertableaux} \ydiagram
[ *(yellow!70) ]{ 1 , 2 } + \ytableausetup {mathmode, boxsize=0.6em,centertableaux} \ydiagram [ *(yellow!70) ]{ 1 } + \ytableausetup
{mathmode, boxsize=0.6em,centertableaux} \ydiagram [ *(yellow!70) ]{ 1 }
+ \U
\end{align*}
\begin{align*}
\vf\left(\ytableausetup {mathmode, boxsize=0.6em,centertableaux} \ydiagram [
*(yellow!70) ]{ 2, 2 , 5 , 5 }
\right)&=\Ct{5}+\Ct{4}+\Ct{3}\cdot\Ct{2}+2\Ct{3}+2\Ct{2}+1&\\
&=42+14+10+10+4+1=81&
\end{align*}

\subsection{Desarrollo de [7,7,2,2]}\label{it8}

\

\

De forma similar al anterior hacen falta 8 iteraciones para que el diagrama queda reducido y transformado a un polinomio en números de Catalan como sigue:

\begin{align*}
\H\left(\ytableausetup {mathmode, boxsize= 0.6 em,centertableaux} \ydiagram[*(red)]
  {1+1,2+0,3+4,4+3}* [
*(yellow!70) ]{ 2,2,7,7} 
\right)
=&3\Ct{2}\Ct{3} + 6\Ct{2} + 5\Ct{3}^2  3\Ct{4}+ \Ct{5} +4&\\
=&155&
\end{align*}

\subsection{Desarrollo de $\D_{8,14}$ [12,10,8,7,5,3,1]}\label{it18}

\

\

Análogamente a los ejemplos anteriores, en este necesitamos18 iteraciones y el diagrama queda reducido y transformado a un polinomio en números de Catalan como sigue:

\begin{align*}
\H\left(\ytableausetup {mathmode, boxsize= 0.6 em,centertableaux} \ydiagram[*(red)]
  {1+0,2+1,3+2,4+3,5+3,6+4,7+5}* [
*(yellow!70) ]{ 1 , 3 , 5 , 7 , 8 , 10 , 12 } 
\right)
=&3\Ct{2}^4 + 11\Ct{2}^3 + 32\Ct{2}^2\Ct{3} + 10\Ct{2}\Ct{3}^2 + 11\Ct{2}^2\Ct{4}+&\\
&+ 3\Ct{2}^2 + 22\Ct{2}\Ct{3} + 17\Ct{3}^2 + 31\Ct{2}\Ct{4} + 22\Ct{3}\Ct{4} +&\\
&+3\Ct{4}^2 + 20\Ct{2}\Ct{5} + 5\Ct{3}\Ct{5} + 5\Ct{2}\Ct{6} + \Ct{3} + 5\Ct{4} +&\\
&+10\Ct{5} + 10\Ct{6} + 5\Ct{7} + \Ct{8}&\\
=&14985&
\end{align*}

\section{Árbol de Kréwéras}
El árbol de Kréwéras relaciona en forma inclusiva los caminos de Dyck a travez de sus diagramas de Ferrers. Comenzando con el diagrama de Christoffel y descendiendo pasando por todos los diagramas de Dyck hasta el diagrama vacío. Por ejemplo el árbol de Kréwéras para $\D_{4,6}$ se ve en la figura \ref{D46}.

\begin{figure}[h]
\begin{align*}
\begin{tikzpicture}[xscale=0.7,yscale=0.5, every node/.style={scale=0.6}] \draw
[color=red!50,line width=0.500000000000000pt,-] (-1,0) -- (-1, 2 );
\draw [color=red!50,line width=0.500000000000000pt,-] ( -1 , 2 ) -- ( -2
, 4 ); \draw [color=red!50,line width=0.500000000000000pt,-] ( -1 , 2 )
-- ( 0 , 4 ); \draw [color=red!50,line width=0.500000000000000pt,-] ( -2
, 4 ) -- ( -3 , 6 ); \draw [color=red!50,line
width=0.500000000000000pt,-] ( -2 , 4 ) -- ( -1 , 6 ); \draw
[color=red!50,line width=0.500000000000000pt,-] ( 0 , 4 ) -- ( -1 , 6 );
\draw [color=red!50,line width=0.500000000000000pt,-] ( 0 , 4 ) -- ( 1 ,
6 ); \draw [color=red!50,line width=0.500000000000000pt,-] ( -3 , 6 ) --
( -4 , 8 ); \draw [color=red!50,line width=0.500000000000000pt,-] ( -3 ,
6 ) -- ( -2 , 8 ); \draw [color=red!50,line width=0.500000000000000pt,-]
( -1 , 6 ) -- ( -2 , 8 ); \draw [color=red!50,line
width=0.500000000000000pt,-] ( -1 , 6 ) -- ( 0 , 8 ); \draw
[color=red!50,line width=0.500000000000000pt,-] ( -1 , 6 ) -- ( 2 , 8 );
\draw [color=red!50,line width=0.500000000000000pt,-] ( 1 , 6 ) -- ( 2 ,
8 ); \draw [color=red!50,line width=0.500000000000000pt,-] ( -4 , 8 ) --
( -4 , 10 ); \draw [color=red!50,line width=0.500000000000000pt,-] ( -2
, 8 ) -- ( -4 , 10 ); \draw [color=red!50,line
width=0.500000000000000pt,-] ( -2 , 8 ) -- ( -2 , 10 ); \draw
[color=red!50,line width=0.500000000000000pt,-] ( -2 , 8 ) -- ( 0 , 10
); \draw [color=red!50,line width=0.500000000000000pt,-] ( 0 , 8 ) -- (
-2 , 10 ); \draw [color=red!50,line width=0.500000000000000pt,-] ( 0 , 8
) -- ( 2 , 10 ); \draw [color=red!50,line width=0.500000000000000pt,-] (
2 , 8 ) -- ( 0 , 10 ); \draw [color=red!50,line
width=0.500000000000000pt,-] ( 2 , 8 ) -- ( 2 , 10 ); \draw
[color=red!50,line width=0.500000000000000pt,-] ( -4 , 10 ) -- ( -4 , 12
); \draw [color=red!50,line width=0.500000000000000pt,-] ( -4 , 10 ) --
( -2 , 12 ); \draw [color=red!50,line width=0.500000000000000pt,-] ( -2
, 10 ) -- ( -4 , 12 ); \draw [color=red!50,line
width=0.500000000000000pt,-] ( -2 , 10 ) -- ( 0 , 12 ); \draw
[color=red!50,line width=0.500000000000000pt,-] ( -2 , 10 ) -- ( 2 , 12
); \draw [color=red!50,line width=0.500000000000000pt,-] ( 0 , 10 ) -- (
-2 , 12 ); \draw [color=red!50,line width=0.500000000000000pt,-] ( 0 ,
10 ) -- ( 2 , 12 ); \draw [color=red!50,line
width=0.500000000000000pt,-] ( 2 , 10 ) -- ( 2 , 12 ); \draw
[color=red!50,line width=0.500000000000000pt,-] ( -4 , 12 ) -- ( -3 , 14
); \draw [color=red!50,line width=0.500000000000000pt,-] ( -4 , 12 ) --
( -1 , 14 ); \draw [color=red!50,line width=0.500000000000000pt,-] ( -2
, 12 ) -- ( -1 , 14 ); \draw [color=red!50,line
width=0.500000000000000pt,-] ( 0 , 12 ) -- ( -3 , 14 ); \draw
[color=red!50,line width=0.500000000000000pt,-] ( 0 , 12 ) -- ( 1 , 14
); \draw [color=red!50,line width=0.500000000000000pt,-] ( 2 , 12 ) -- (
-1 , 14 ); \draw [color=red!50,line width=0.500000000000000pt,-] ( 2 ,
12 ) -- ( 1 , 14 ); \draw [color=red!50,line
width=0.500000000000000pt,-] ( -3 , 14 ) -- ( -1 , 16 ); \draw
[color=red!50,line width=0.500000000000000pt,-] ( -1 , 14 ) -- ( -1 , 16
); \draw [color=red!50,line width=0.500000000000000pt,-] ( 1 , 14 ) -- (
-1 , 16 ); \draw[color= red ]( -1 , 0 ) node{ $\U$ }; \draw [color= red!50 ]( -1
, 2 ) node{$ \bullet$}; \draw ( -1 , 2 ) node{ \ytableausetup {mathmode,
boxsize= 1 em,centertableaux} \ydiagram [ *(red) ]{ 1 } }; \draw [color=
red!50 ]( -2 , 4 ) node{$ \bullet$}; \draw ( -2 , 4 ) node{
\ytableausetup {mathmode, boxsize= 1 em,centertableaux} \ydiagram [
*(red) ]{ 2 } }; \draw [color= red!50 ]( 0 , 4 ) node{$ \bullet$}; \draw
( 0 , 4 ) node{ \ytableausetup {mathmode, boxsize= 1 em,centertableaux}
\ydiagram [ *(red) ]{ 1 , 1 } }; \draw [color= red!50 ]( -3 , 6 ) node{$
\bullet$}; \draw ( -3 , 6 ) node{ \ytableausetup {mathmode, boxsize= 1
em,centertableaux} \ydiagram [ *(yellow!70) ]{ 3 } }; \draw [color=
red!50 ]( -1 , 6 ) node{$ \bullet$}; \draw ( -1 , 6 ) node{
\ytableausetup {mathmode, boxsize= 1 em,centertableaux} \ydiagram [
*(red) ]{ 1 , 2 } }; \draw [color= red!50 ]( 1 , 6 ) node{$ \bullet$};
\draw ( 1 , 6 ) node{ \ytableausetup {mathmode, boxsize= 1
em,centertableaux} \ydiagram [ *(red) ]{ 1 , 1 , 1 } }; \draw [color=
red!50 ]( -4 , 8 ) node{$ \bullet$}; \draw ( -4 , 8 ) node{
\ytableausetup {mathmode, boxsize= 1 em,centertableaux} \ydiagram [
*(yellow!70) ]{ 4 } }; \draw [color= red!50 ]( -2 , 8 ) node{$
\bullet$}; \draw ( -2 , 8 ) node{ \ytableausetup {mathmode, boxsize= 1
em,centertableaux} \ydiagram [ *(yellow!70) ]{ 1 , 3 } }; \draw [color=
red!50 ]( 0 , 8 ) node{$ \bullet$}; \draw ( 0 , 8 ) node{ \ytableausetup
{mathmode, boxsize= 1 em,centertableaux} \ydiagram [ *(red) ]{ 2 , 2 }
}; \draw [color= red!50 ]( 2 , 8 ) node{$ \bullet$}; \draw ( 2 , 8 )
node{ \ytableausetup {mathmode, boxsize= 1 em,centertableaux} \ydiagram
[ *(red) ]{ 1 , 1 , 2 } }; \draw [color= red!50 ]( -4 , 10 ) node{$
\bullet$}; \draw ( -4 , 10 ) node{ \ytableausetup {mathmode, boxsize= 1
em,centertableaux} \ydiagram [ *(yellow!70) ]{ 1 , 4 } }; \draw [color=
red!50 ]( -2 , 10 ) node{$ \bullet$}; \draw ( -2 , 10 ) node{
\ytableausetup {mathmode, boxsize= 1 em,centertableaux} \ydiagram [
*(yellow!70) ]{ 2 , 3 } }; \draw [color= red!50 ]( 0 , 10 ) node{$
\bullet$}; \draw ( 0 , 10 ) node{ \ytableausetup {mathmode, boxsize= 1
em,centertableaux} \ydiagram [ *(yellow!70) ]{ 1 , 1 , 3 } }; \draw
[color= red!50 ]( 2 , 10 ) node{$ \bullet$}; \draw ( 2 , 10 ) node{
\ytableausetup {mathmode, boxsize= 1 em,centertableaux} \ydiagram [
*(red) ]{ 1 , 2 , 2 } }; \draw [color= red!50 ]( -4 , 12 ) node{$
\bullet$}; \draw ( -4 , 12 ) node{ \ytableausetup {mathmode, boxsize= 1
em,centertableaux} \ydiagram [ *(yellow!70) ]{ 2 , 4 } }; \draw [color=
red!50 ]( -2 , 12 ) node{$ \bullet$}; \draw ( -2 , 12 ) node{
\ytableausetup {mathmode, boxsize= 1 em,centertableaux} \ydiagram [
*(yellow!70) ]{ 1 , 1 , 4 } }; \draw [color= red!50 ]( 0 , 12 ) node{$
\bullet$}; \draw ( 0 , 12 ) node{ \ytableausetup {mathmode, boxsize= 1
em,centertableaux} \ydiagram [ *(yellow!70) ]{ 3 , 3 } }; \draw [color=
red!50 ]( 2 , 12 ) node{$ \bullet$}; \draw ( 2 , 12 ) node{
\ytableausetup {mathmode, boxsize= 1 em,centertableaux} \ydiagram [
*(yellow!70) ]{ 1 , 2 , 3 } }; \draw [color= red!50 ]( -3 , 14 ) node{$
\bullet$}; \draw ( -3 , 14 ) node{ \ytableausetup {mathmode, boxsize= 1
em,centertableaux} \ydiagram [ *(yellow!70) ]{ 3 , 4 } }; \draw [color=
red!50 ]( -1 , 14 ) node{$ \bullet$}; \draw ( -1 , 14 ) node{
\ytableausetup {mathmode, boxsize= 1 em,centertableaux} \ydiagram [
*(yellow!70) ]{ 1 , 2 , 4 } }; \draw [color= red!50 ]( 1 , 14 ) node{$
\bullet$}; \draw ( 1 , 14 ) node{ \ytableausetup {mathmode, boxsize= 1
em,centertableaux} \ydiagram [ *(yellow!70) ]{ 1 , 3 , 3 } }; \draw
[color= red!50 ]( -1 , 16 ) node{$ \bullet$}; \draw ( -1 , 16 ) node{
\ytableausetup {mathmode, boxsize= 1em,centertableaux} \ydiagram [
*(yellow!70) ]{ 1 , 3 , 4 } }; \end{tikzpicture}
\end{align*}
\caption{Árbol de Kréwéras de $\D_{4,6}$.}\label{D46}
\end{figure}
Hemos coloreado en rojo todos los diagramas relacionados con el diagrama $[2,2,1]$ (ver ejemplo \ref{e221}). El numero de diagramas en la rama es exactamente igual al numero de caminos de Dyck, $\H\left(\ytableausetup {mathmode, boxsize= 0.4em,centertableaux} \ydiagram [*(yellow!70) ]{ 1,2 , 2 }\right)=9$. Debido al hecho que el numero de elementos de una rama del árbol es igual al numero de caminos de Dyck contenidos en el diagrama podemos contar dichos elementos a partir de nuestro método de descomposición de diagramas.\\
\section{Conclusión}
Como hemos visto el método de descomposición de diagramas resuelve de forma simple y recursiva la obtención del numero de caminos de Dyck en los casos mas generales. La desventajas es no conocer aun la forma mas efectiva de reducir las iteraciones ya que para algunos casos el numero de iteraciones puede ser muy elevado. Esto puede verse en los ejemplos  \ref{it8b}, \ref{it8} y \ref{it18}. Una opción es el método alternativo que puede disminuir la cantidad de iteraciones en casos particulares pero la descripción es como diferencia y no como suma (ver ejemplos  \ref{311a}). Hay muchas derivaciones de los métodos como por ejemplo encontrar el polinomio q-analógico. Todos los algoritmos en SAGE están disponible en \url{http://thales.math.uqam.ca/~jeblazek/Sage_Combinatory.html}

\section*{Agradecimientos}
Deseo agradecer a mi director François Bergeron, a Sre\v{c}ko Brlek, a Hector Blandin Noguera y a Yannic Vargas Lozada, por sus consejos y por su apoyo durante la preparación de este documento. 
\newpage

\bibliographystyle{splncs03.bst}
\bibliography{Bibliog}

\begin{thebibliography}{10}
\providecommand{\url}[1]{\texttt{#1}}
\providecommand{\urlprefix}{URL }

\bibitem{BH09}
Bakhtin, Y., Heitsch, C.E.: Large deviations for random trees and the branching
  of {RNA} secondary structures. Bull. Math. Biol.  71(1),  84--106 (2009),
  \url{http://dx.doi.org/10.1007/s11538-008-9353-y}

\bibitem{BLL}
Bergeron, F., Labelle, G., Leroux, P., Readdy, M.: Combinatorial species and
  tree-like structures. Encyclopedia of mathematics and its applications,
  Cambridge University Press, Cambridge, England (1998),
  \url{http://opac.inria.fr/record=b1098109}

\bibitem{Biz54}
Bizley, M.T.L.: Derivation of a new formula for the number of minimal lattice
  paths from $(0,0)$ to $(km,kn) \cdots$. JIA 80 pp. 55--62 (1954)

\bibitem{BJE14}
Blazek, J.E.: Conbinatoire de $\N$-modules {C}atalan, Master Thesis.
  D\'epartement de Math\'ematique, UQAM (2015),
  \url{http://thales.math.uqam.ca/~jeblazek/Articles/Memoire2014A.pdf}

\bibitem{Chap09}
Chapman, R.J., Chow, T.Y., Khetan, A., Moulton, D.P., Waters, R.J.: Simple
  formulas for lattice paths avoiding certain periodic staircase boundaries. J.
  Combin. Theory Ser. A  116(1),  205--214 (2009),
  \url{http://dx.doi.org/10.1016/j.jcta.2008.05.002}

\bibitem{Duchon}
Duchon, P.: On the enumeration and generation of generalized dyck words.
  Discrete Mathematics  225(1-3),  121--135 (2000)

\bibitem{Eilenberg}
Eilenberg, S.: Automata, Languages, and Machines. Academic Press, Inc.,
  Orlando, FL, USA (1976)

\bibitem{Fu13}
Fukukawa, Y.: Counting generalized {D}yck paths. arXiv:1304.5595v1 [math.CO]
  (April 2013)

\bibitem{GMV14}
Gorky, E., Mazin, M., Vazirani, M.: Affine permutations and rational slope
  parking functions. arXiv:1403.0303v1 [math.CO]  (Mars 2014)

\bibitem{Serrano03}
Goulden, I.P., Serrano, L.G.: Maintaining the spirit of the reflection
  principle when the boundary has arbitrary integer slope. J. Combin. Theory
  Ser. A  104(2),  317--326 (2003),
  \url{http://dx.doi.org/10.1016/j.jcta.2003.09.004}

\bibitem{KOS09}
Koshy, T.: Catalan numbers with applications. Oxford University Press, Inc.
  (2009)

\bibitem{Krew65}
Kréwéras, G.: Sur une classe de problèmes de dénombrement liés au treillis
  des partitions des entiers. No.~6, Institut de statistique des universités
  de Paris, Paris, France (1965)

\bibitem{LY90}
Labelle, J., Yeh, Y.: Generalized dyck paths. Discrete Mathematics  82(1),
  1--6 (1990)

\bibitem{Lothaire1}
Lothaire, M.: Applied Combinatorics on Words. Cambridge University Press (2005)

\bibitem{MR13}
Melan{\c c}on, G., Reutenauer, C.: On a class of {L}yndon words extending
  {C}hristoffel words and related to a multidimensional continued fraction
  algorithm. Journal of Integer Sequences  16(Article 13.9.7) (2013)

\bibitem{SW94}
Schmitt, W.R., Waterman, M.S.: Linear trees and {RNA} secondary structure.
  Discrete Appl. Math.  51(3),  317--323 (1994),
  \url{http://dx.doi.org/10.1016/0166-218X(92)00038-N}

\bibitem{Sim00}
Simion, R.: Noncrossing partitions. Discrete Math.  217(1-3),  367--409 (2000),
  \url{http://dx.doi.org/10.1016/S0012-365X(99)00273-3}, formal power series
  and algebraic combinatorics (Vienna, 1997)

\bibitem{Stanley}
Stanley, R.P.: Enumerative Combinatorics: Volume 1. Cambridge University Press,
  New York, NY, USA, 2nd edn. (2011)

\end{thebibliography}

\end{document}